\documentstyle[12pt]{amsart} \mathsurround=2pt
\newtheorem{thm}{Theorem}[section]

\newtheorem{lm}[thm]{Lemma}

\newcommand{\der}{\partial}

\begin{document}

\title{Codimension growth  of a variety of Novikov algebras}
\author{A.S. Dzhumadil'daev}
\address{Institute of Mathematics, Pushkin street 125, Almaty, KAZAKHSTAN }
\email{askar56@@hotmail.com}

\subjclass{Primary 16R10,  17A50, 17A30, 17D25, 17C50}

\keywords{Novikov algebra, codimensions sequence, polynomial identities, growth of a variety}

\maketitle
\begin{abstract} An algebra with identities $a\circ(b\circ c-c\circ b)=(a\circ b)\circ c-(a\circ c)\circ b$
and $a\circ(b\circ c)=b\circ(a\circ c)$ is called Novikov. We construct free Novikov base in terms of Young diagrams. We show that codimensions exponent for 
a variety of Novikov algebras exists and is equal $4$.
\end{abstract}

\section{Introduction}

A variety of algebras is a class of algebras satisfying some polynomial identities. 
One of important parameters of varieties is so-called codimension growth. 
If $\mathcal V$ is a variety and $N_n({\mathcal V})$ is a polylinear part of its free algebra  generated by $n$ elements, then $c_n({\mathcal V})=dim\,N_n(\mathcal V)$ is called $n$-th codimension of $\mathcal V$. Codimension growth is defined by  a sequence of codimensions $c_1,c_2,c_3,\ldots.$ {\it Codimension exponent} is defined as 
$$Exp({\mathcal V})=\lim c_n({\mathcal V})^{1/n}.$$
Appear natural questions whether this exponent exists and  whether it is an integer. In associative case these questions are well studied. It was proved that $Exp({\mathcal V})$ exists and it is an integer for any proper variety of associative algebras 
(\cite{Zaicev}).  Constructions of free bases for Lie algebras are well known (about Hall-Lyndon- Shrishov bases see, for example \cite{FreeLie}).

In our paper we consider class of non-associative algebras. 
An algebra $A=(A,\circ)$ is called {\it Novikov} (\cite{Novikov}, \cite{Osborn}, \cite{Gelfand}), if it satisfies the following identities 
$$(a,b,c)=(a,c,b)$$
$$a\circ(b\circ c)=b\circ(a\circ c),$$
for any $a,b,c\in A.$ Here 
$$(a,b,c)=a\circ(b\circ c)-(a\circ b)\circ c$$
is an associator. Novikov algebras are Lie-admissible. 

{\bf Example.} $A={\bf C}[x]$ under multiplication $a\circ b=\der(a)b$ is Novikov. 

Let $Nov$ be a variety of Novikov algebras and $N_n$ be a polylinear part of free Novikov algebra generated by $n$ elements. Let    
$$ Exp(Nov){=}\lim_{n\rightarrow \infty}
(\dim N_n)^{1/n}$$
be codimensions growth of Novikov variety. 

In \cite{Dzhuma} was proved that Novikov operad is not Koszul. 
Now we give construction of free Novikov base in terms of Young diagrams and use this base to calculation of  generating function and codimension growth.  Our  main result is 

\begin{thm} \label{polypart} Codimensions sequence for Novikov varety is given by 
$$\dim N_n={2n-2\choose n-1}.$$
Codimension exponent of Novikov variety exists and  $$Exp(Nov)=4.$$ Generating function of codimensions sequence $\sum_{i\ge 0} N_ix^i$
is equal to $x(1-4 x)^{-1/2}.$
\end{thm}

\section{\label{basic} Free base for  Novikov
algebras}

In \cite{Dzh1} are given constructions of a base of free Novikov
algebra in terms of $r$-elements and in terms of rooted trees. In
this section we give  construction of free  base in
terms of Young diagrams.

Recall that Young diagram is a set of boxes (we denote them by
bullets) with non-increasing numbers of boxes in each row. Rows
and columns  are numerated from top to bottom and from left to
right. Let $k$ be the  number of rows and $r_i$ be the number of
boxes in $i$-th row. The total number of boxes, $r_1+\cdots+r_k,$
is called {\it degree} of Young diagram.

To construct Novikov diagram, we need to complement Young diagram
by one box, we call it as "a nose". Namely, we need to add the first row by one more box,
$$\begin{array}{ccccc} \bullet&\cdots&\bullet&\bullet&\bullet\\
\bullet&\cdots&\bullet&\bullet&\\ \vdots& \cdots&\vdots &\vdots&\\
\bullet&\cdots&\bullet&&\\
\end{array}
\mapsto
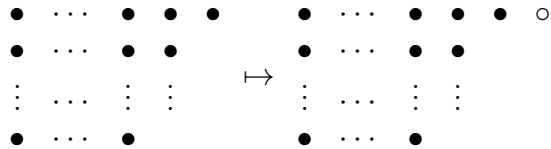
\begin{array}{cccccc}
\bullet&\cdots&\bullet&\bullet&\bullet&\circ\\
\bullet&\cdots&\bullet&\bullet&\\ \vdots&
\cdots&\vdots&\vdots&&\\\bullet&\cdots&\bullet&&&\\
\end{array}
$$ The number of boxes in Novikov diagram is called its {\it degree}. So,
difference between degrees of Novikov diagram and corresponding
Young diagram is equal to $1.$

Let us given an alphabet (ordered set) $\Omega.$ To construct
Novikov tableau on $\Omega$ we need to feel  Novikov diagrams by
elements of $\Omega.$ Denote by $a_{i,j}$ an element of $\Omega$
in the box $(i,j),$ that is a cross of $i$-th row by $j$-th
column. The feeling rule is the following
\begin{itemize}
\item $a_{i,1}\ge a_{i+1,1},$ if $r_i=r_{i+1},
i=1,2,\ldots,k-1.$
\item the sequence $a_{k,2}\cdots a_{k,r_k}
a_{k-1,2}\cdots a_{k-1,r_{k-1}}\cdots a_{1,2}\cdots
a_{1,r_1}a_{1,r_1+1}$  is non-decreasing.
\end{itemize}
In particular, all boxes beginning from the second place in each
row are labeled by non-decreasing elements of alphabet. Denote by
$R_n$ a set of Novikov tableaux labeled by $\Omega$ with
$n=r_k+\cdots +r_1+1$ boxes.

Let $F(\Omega)$ be free Novikov algebra generated by $\Omega.$ Let
$F_n(\Omega)$ be its subspace generated by basic elements of
degree $n.$ Correspond to any Novikov tableaux
$$\begin{array}{cccccc}
a_{1,1}&\cdots&\cdots&a_{1,r_1-1}&a_{1,r_1}&a_{1,r_1+1}\\
a_{2,1}&\cdots&a_{2,r_2-1}&a_{2,r_2}&\\ \vdots&
\cdots&\vdots&\vdots&&\\ a_{k,1}&\cdots&a_{k,r_k}&&&\\
\end{array}$$
an element $$X= X_k\circ (X_{k-1}\circ( \cdots \circ (X_2\circ
X_1)\cdots )),$$ (right-normed bracketing) where $$X_i=(\cdots
((a_{i,1}\circ a_{i,2})\circ a_{i,3})\cdots \circ
a_{i,r_i-1})\circ a_{i,r_i},\quad 1<i\le k,$$ $$X_1=(\cdots
((a_{1,1}\circ a_{1,2})\circ a_{1,3})\cdots \circ a_{1,r_1})\circ
a_{1,r_1+1}$$ (left-normed bracketing). All base elements of free
Novikov algebra $F(\Omega)$ are obtained by this way. In
particular, $\dim F_n(\Omega)=|R_n|.$

As an example let us construct base elements of polylinear part $N_4$ of free Novikov algebra generated by $4$ elements $a,b,c,d.$

Young diagrams of degree  $3$:

$$\begin{array}{c}\bullet\\  \bullet\\  \bullet\\
\end{array}\qquad
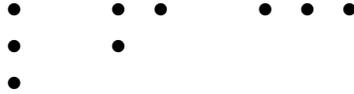
\begin{array}{cc}\bullet&\bullet\\
\bullet&\\&\\
\end{array}
\qquad
\begin{array}{ccc}
\bullet&\bullet&\bullet\\ &&\\ &&\\
\end{array}
$$

{\noindent Novikov diagrams of degree 4:}

$$
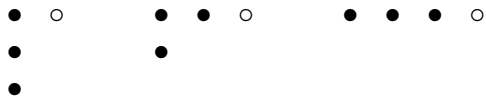
\begin{array}{cc}\bullet&\circ\\\bullet&\\\bullet&\\
\end{array}\qquad
\begin{array}{ccc}\bullet&\bullet&\circ\\
\bullet&&\\ &&\\
\end{array}\qquad
\begin{array}{cccc}
\bullet&\bullet&\bullet&\circ\\ &&&\\&&&\\
\end{array}
$$

{\noindent Novikov tableaux of degree 4 generated by elements $a,b,c,d$:}

$$
\begin{array}{cc}c&d\\b&\\a&\\ \end{array}\qquad
\begin{array}{cc}d&c\\b&\\a&\\ \end{array}\qquad
\begin{array}{cc}d&b\\c&\\a&\\ \end{array}\qquad
\begin{array}{cc}d&a\\c&\\b&\\ \end{array}
$$

\bigskip
\bigskip

$$
\begin{array}{ccc}b&c&d\\
a&&\\ &&\\
\end{array}\qquad
\begin{array}{ccc}c&b&d\\
a&&\\ &&\\
\end{array}\qquad
\begin{array}{ccc}d&b&c\\
a&&\\ &&\\
\end{array}
$$ $$
\begin{array}{ccc}a&c&d\\
b&&\\ &&\\
\end{array}\qquad
\begin{array}{ccc}c&a&d\\
b&&\\ &&\\
\end{array}\qquad
\begin{array}{ccc}d&a&c\\
b&&\\ &&\\
\end{array}
$$ $$
\begin{array}{ccc}a&b&d\\
c&&\\ &&\\
\end{array}\qquad
\begin{array}{ccc}b&a&d\\
c&&\\ &&\\
\end{array}\qquad
\begin{array}{ccc}d&a&b\\
c&&\\ &&\\
\end{array}
$$ $$
\begin{array}{ccc}a&b&c\\
d&&\\ &&\\
\end{array}\qquad
\begin{array}{ccc}b&a&c\\
d&&\\ &&\\
\end{array}\qquad
\begin{array}{ccc}c&a&b\\
d&&\\ &&\\
\end{array}
$$

\bigskip

$$\begin{array}{cccc} a&b&c&d\\ &&&\\&&&\\
\end{array}\qquad
\begin{array}{cccc} b&a&c&d\\ &&&\\&&&\\
\end{array}\qquad
\begin{array}{cccc} c&a&b&d\\ &&&\\&&&\\
\end{array}\qquad
\begin{array}{cccc} d&a&b&c\\ &&&\\&&&\\
\end{array}\qquad
$$

{\noindent So,  polylinear part of free Novikov algebra in degree 4 is $20$-dimensional and 
the following elements form base } 

$$a\circ(b\circ (c\circ d)), a\circ(b\circ (d\circ c)),
a\circ(c\circ (d\circ b)), b\circ(c\circ (d\circ a)),$$

$$a\circ ((b\circ c)\circ d), a\circ ((c\circ b)\circ d),
a\circ((d\circ b)\circ c),$$ $$b\circ ((a\circ c)\circ d), b\circ
((c\circ a)\circ d), b\circ((d\circ a)\circ c),$$ $$c\circ
((a\circ b)\circ d), c\circ ((b\circ a)\circ d), c\circ((d\circ
a)\circ b),$$ $$d\circ ((a\circ b)\circ c), d\circ ((b\circ
a)\circ c), d\circ((c\circ a)\circ b),$$

$$((a\circ b)\circ c)\circ d, ((b\circ a)\circ c)\circ d, ((c\circ
a)\circ b)\circ d, ((d\circ a)\circ b)\circ c.$$

\section{Codimensions growth of Novikov variety}

Let $\lambda=1^{m_1(\lambda)}2^{m_2(\lambda)}\cdots $ be a partition
of $n-1,$ i.e., $$|\lambda|=\sum_{i\ge 1} i\,m_i(\lambda)=n-1.$$
Let $$m(\lambda)=\sum_{i\ge 1}m_i(\lambda).$$ For
$m,m_1,m_2,\ldots,m_n,$ such that $m=m_1+m_2+\cdots+m_n$ let
$${m\choose m_1,m_2,\ldots ,m_n}=\frac{m!}{m_1!m_2!\cdots m_n!}$$
be a multinomial coefficient.

\begin{lm}\label{kuralai} $$\sum_{|\lambda|=n-1}{m(\lambda)\choose
m_1(\lambda), m_2(\lambda), \ldots}{n\choose
m(\lambda)}={2(n-1)\choose n-1}$$
\end{lm}

{\bf Proof.} By Vandermonde convolution relation \cite{Riordan},
chapter 1.3, formulae (3a), $${n+p\choose m}=\sum_{s\ge 0}
{n\choose m-s}{p\choose s}.$$

By \cite{Riordan} chapter 4.5, formulae (21), for fixed $n$ and
$m$, takes place the following relation
$$\sum_{m_1,m_2,\ldots}{m\choose m_1,m_2,\ldots ,m_n}={n-1\choose
m-1},$$ where summation is over $m_1,m_2,\ldots,m_n,$ such that
$m=m_1+m_2+\cdots+m_n, n=m_1+2m_2+\cdots+n m_n.$

By these relations,
 $$\sum_{|\lambda|=n-1}{m(\lambda)\choose
m_1(\lambda), m_2(\lambda), \ldots}{n\choose m(\lambda)}$$
\medskip
$$=\sum_{s\ge 1} {n\choose s}\sum_{|\lambda|=n-1,
m(\lambda)=s}{s\choose m_1(\lambda), m_2(\lambda), \ldots}$$
\medskip
$$=\sum_{s\ge 1} {n\choose s} {n-2\choose s-1}=\sum_{s\ge 1}
{n\choose s} {n-2\choose n-s-1}$$ 
\medskip
$$={2(n-1)\choose n-1}.$$

\begin{lm}\label{arman}
$$\lim_{n\rightarrow \infty} {2n-2\choose n-1}^{1/n}=4.$$
\end{lm}

{\bf Proof.} 
For $a\in {\bf Z}$ denote by $a!!$ a product of positive integers
$a, a-2, a-4,$ and so on. For example, $(2n-1)!!$ is a product of
odd numbers between $1$ and $2n-1.$ We have $$(2n-4)!!\le
(2n-3)!!\le (2n-2)!!.$$ Thus $$\frac{2^{n-2}}{n-1}\le
\frac{(2n-3)!!}{(n-1)!}\le 2^{n-1}.$$ Since $${2(n-1)\choose
n-1}=\frac{2(n-1)!}{((n-1)!)^2}=\frac{2^{n-1}(2n-3)!!}{(n-1)!},$$ we have
$$\frac{2^{2n-3}}{n-1}\le {2(n-1)\choose n-1}\le 2^{2(n-1)}.$$
It remains to note that  $$\lim_{n\rightarrow \infty}
(\frac{2^{2n-3}}{(n-1)})^{1/n} =\lim_{n\rightarrow
\infty}(2^{4n-2})^{1/n}=4.$$

\bigskip

{\bf Proof of Theorem \ref{polypart}.}

As we mentioned above any polylinear base element of free Novikov
algebra of degree $n$ corresponds to Young diagram of degree $n-1$.
Suppose that it has all together $m$ rows, namely $m_1$ rows with
$i_1$ boxes, $m_2$ rows with $i_2$ boxes, etc, $m_k$ rows with
$i_k$ boxes, where $i_1>i_2>\cdots
>i_k.$  So, $\sum_{s=1}^k i_sm_s=n-1,$ and such Young diagram looks like
 $$\begin{array}{lll} m_1&\{&{\begin{array}{cccccc}
\bullet&\cdots&\bullet&\bullet&\bullet&\bullet\\ \vdots &\cdots
&\vdots&\vdots&\vdots&\vdots\\
\bullet&\cdots&\bullet&\bullet&\bullet&\bullet
\end{array}}\\
 m_2&\{&
\begin{array}{cccccc}
\bullet&\cdots&\bullet&\bullet&&\\ \vdots&\cdots&\vdots&\vdots&&\\
\bullet&\cdots&\bullet&\bullet&&
\end{array}\\
&&\begin{array}{cccccc} \vdots&\;\cdots&\vdots&&&\\\end{array}\\
m_k&\{&\begin{array}{cccccc} \bullet&\cdots&\bullet&&&\\
\vdots&\cdots&\vdots&&&\\
 \bullet&\cdots&\bullet&&&\\
\end{array}
\end{array}
$$ Set $m_i=0,$ if $i>k.$ The Novikov diagram corresponding to
such Young diagram, filled by $n$ different letters, is uniquiely
defined by its first column. The first column can be choosen in
$${n\choose m_1,m_2,\ldots,m_n}={m\choose m_1,m_2,\ldots,m_n}
{n\choose m}$$ ways. Therefore, by lemma \ref{kuralai} $$\dim
R_n=\sum {m\choose m_1,m_2,\ldots,m_n} {n\choose m}={2(n-1)\choose
n-1}$$ (summation is over $m_1,m_2,\ldots,m_n$ such that
$\sum_{s}i_sm_s=n-1.$)
By Lemma \ref{arman}, $Exp(Nov)=4.$
Other statements of theorem \ref{polypart} is evident.


\begin{thebibliography}{20}

\bibitem{Novikov}  A.A. Balinskii,   S.P. Novikov, {\em Poisson bracket
of hamiltonian type, Frobenius algebras and Lie algebras,} Dokladu
AN SSSR {\bf 1985}, {\it 283} (5), 1036--1039.

\bibitem{Dzh1}  A.S. Dzhumadil'daev, C. Lofwall, 
{\em Trees, free right-symmetric algebras, free Novikov algebras
and identities}, Homology, Homotopy and Appl., {\bf 4}\{2002),
No.2(1), 165-190.

\bibitem{Dzhuma} A.S. Dzhumadil'daev, 
{\em  Algebras with skew-symmetric identity of degree $3$}, 
Sovremennaya Matematika i ee Prilozheniya (VINITI) 60 (2008), 13--31. 
(English translation to appear in J. Math. Sci. )

\bibitem{Gelfand}  I.M. Gelfand,  I.Ya. Dorfman,   {\em
Hamiltonian operators and related algebraic structures,} Funct.
Anal. Prilozhen. {\bf 1979}, {\it 13} (4), 13--30= engl.transl.
Funct. Anal. Appl. {\bf 1979}, {\it 13}, 248--262.

\bibitem{Zaicev}  A. Giambruno,  M. Zaicev, {\em Minimal varieties of algebras of exponential growth,}
Elect. Res. Ann.  AMS, {\bf 6}(2000), 40-44.

\bibitem{Osborn}  J.M.  Osborn,   {\em Infinite-dimensional Novikov algebras of
characteristic $0$,} J. Algebra {\bf 1994},  {\it 167}, 146--167.

\bibitem{FreeLie}  C. Reutenauer, {\em Free Lie algebras,}  Clarendon Press, Oxford, 1993.

\bibitem{Riordan}  J. Riordan , {\em Combinatorial identities},  J.Wiley\&Sons, Inc., New York-London-Sidney, 1968.

\end{thebibliography}
\end{document}